\def\obs#1{{\bf (*** #1 ***)} }

 \def\obs#1{}     % Remova esta linha para rodar a versao 1

\documentclass[twoside,letterpaper,draft,11pt]{amsart}

\usepackage{amsmath,amsthm,latexsym,xspace,amscd,amssymb,xypic, amscd}               % AmSLaTeX
\usepackage{color}
\usepackage{fullpage}

\thanks{Primary 54H15; Secondary 54H05, 54E50, 54E35.}

\newtheorem{teo}{Theorem}[section]
\newtheorem{defi}[teo]{Definition}
\newtheorem{lema}[teo]{Lemma}

\newtheorem{prop}[teo]{Proposition}

\newtheorem{exe}[teo]{Example}
\newtheorem{rem}[teo]{Remark}

\newcommand{\X}{{\mathbb X}}
\newcommand{\Y}{{\mathbb Y}}
\newcommand{\W}{{\mathbb W}}
\newcommand{\UU}{{\mathbb U}}

\newcommand{\Z}{{\mathbb Z}}

\newcommand{\te}{\theta}

\newcommand{\ti}{{\times}}

\newcommand{\ex}{{\exists}}

\newcommand{\m}{{}^{-1}}
\newcommand{\mt}{\mapsto}

\def\ndv{\ {\mid \kern -0.7 em {\scriptstyle \not}} \ \ }
\def\nd{\ {\mid \kern -0.4 em {\scriptstyle \not}} \ \ }

\newcommand{\N}{{\mathbb N}}

\numberwithin{equation}{section}
\title[Polish globalization of Polish group partial actions]{Polish globalization of Polish group partial actions}

\begin{document}
\author[H.\ Pinedo ]{H. Pinedo }
\address{Escuela de Matem\'aticas, Universidad Industrial de Santander, Cra. 27 Calle 9 UIS Edificio 45\\  Bucaramanga, Colombia}\email{hpinedot@uis.edu.co }
\author[C.\ Uzcategui ]{C. Uzc\'ategui}
\address{Escuela de Matem\'aticas, Universidad Industrial de Santander, Cra. 27 Calle 9 UIS Edificio 45\\  Bucaramanga, Colombia}\email{cuzcatea@uis.edu.co}
\date{\today}

\begin{abstract}  Let  $\X$  be  a separable metrizable  space. We establish a criteria for the existence of a metrizable  globalization for a given continuous partial action of a separable metrizable group $G$  on $\X.$ If $G$ and $\X$ are Polish spaces, we show that the globalization is also a Polish space. We also  show the existence of an universal globalization for partial actions of  Polish groups.
\end{abstract}
\maketitle

\section{Introduction}

Given an action $a:G\times \Y\rightarrow \Y$ of a group $G$ over a set $\Y$ and an invariant  subset $\X$ of $\Y$ (i.e. $a(g,x)\in \X,$ for all $x\in \X$ and $g\in G$),  the restriction of $a$ to $G\times\X$ is an action of $G$ over $\X$. However, if $\X$ is not invariant,  we get what is called a {\em partial action} on $\X$: a collection of partial maps $\{m_g\}_{g\in G}$ on  $\X$ satisfying $m_1={\rm id}_\X$ and $m_g\circ m_h\subseteq m_{gh},$ for all $g,h\in G$.
Partial actions were introduced by R. Exel \cite{E-1,E1}   as an efficient tool for  studying  $C^*$-algebras generated by partial isometries in a Hilbert space, this lead to the  characterizations of various important classes of
 $C^*$-algebras as crossed products by partial actions. For instance,  the Bunce-Deddens and the Bunce-Deddens-Toeplitz algebras \cite{E-2}, the approximately finite
dimensional algebras \cite{E-3}, the Toeplitz algebras of quasi-ordered groups,  as
the Cuntz-Krieger algebras \cite{E-4,QR}, as well as Leavitt-Path algebras \cite{GR, GRO} (see also \cite{D2}).

A natural question is whether  a partial action  can be realized as restrictions  of a corresponding
collection of total maps on some superspace. In the topological context, this problem was
studied  by Abadie \cite{AB} and independently  by J. Kellendonk and M. Lawson  \cite{KL}.  They showed that for any continuous partial action $m$ of a topological group $G$ on a topological space $\X$, there is a topological space $\Y$ and a continuous action $a$ of $G$ on  $\Y$ such that $\X$ is a subspace of $\Y$ and $m$ is the restriction of $a$ to $\X$.  Such a space $\Y$ is called a {\em globalization} of $\X$. They also show that there is a minimal globalization $\X_G$ called the {\em enveloping space} of $\X$.

In general $\X_G$  is not even Hausdorff. In \cite{EG} the authors  shows that the enveloping space of a partial action of a countable  group on compact metric spaces is Hausdorf iff the domain of each $m_g$ is clopen for all  $g\in G$. In \cite{MS}  an approach based on the notion of confluent partial actions is presented. Likewise in \cite{St2}, it was proven that
 the globalization of a partial action on a connected $2$-complex may result in a complex which is not connected, in particular, the enveloping space of $\X$
does not  inherit its structure. So a natural problem is to know which properties of $\X$ are preserved by $\X_G.$

One of the purposes of this article is to study the metrizability of  $\X_G$. We will show that under some  conditions (essentially one that guarantees that $\X_G$ is Hausdorff)  $\X_G$ is metrizable when $G$ and $\X$ are separable metrizable and  is Polish when $\X$ and $G$ are Polish (see Theorem \ref{xgpol}).  This result was known for a countable discrete group  acting partially  on the Cantor space \cite{EG}. Our proof uses a generalization to the context of partial actions of a theorem of Mackey-Hjorth about extensions of Polish group actions (see Theorem \ref{regcuo}) and the fact that $\X_G$ is equal to a quotient space by an orbit equivalence relation (see  Theorem \ref{mtomg}).

Actions of Polish groups have received a lot of attention  in recent  years, because of  its connection with many areas of mathematics (see \cite{KEB,GA} and the references therein).  A particularly interesting feature of Polish group actions is  the existence of a universal one. Becker and Kechris  \cite{KEB} have shown that for every Polish group $G$ there is Polish $G$-space $\UU_G$  such that for every Polish $G$-space $\X,$  there is a (necessarily invariant) set $\Y\subseteq \UU_G$ such that $\X$ and $\Y$ are Borel isomorphic as $G$-spaces.  Since, under suitable conditions, $\X_G$ is Polish when $\X$ and $G$ are Polish,  we conclude  that $\UU_G$  is a universal globalization for  partial actions of $G$ on Polish spaces. In Section \S\ref{universal} we will explore these ideas in the context of partial actions of a countable discrete group.

\section{Terminology}
\label{notions}

We start by recalling some basic notions on descriptive set theory, which can be found in \cite{KE}.  A {\it Polish} space (group) is  a  separable, completely metrizable topological space (group).
A subset of a Polish space is {\it Borel}, if it belongs to the $\sigma$-algebra generated by the open sets.
A function $f\colon \Y\to \X$ is Borel measurable if $f\m(U)$ is Borel in $\Y,$ for any open set $U$ in $\X.$  If $f$ is bijective, Borel measurable and  $f\m$ is also Borel measurable, then $f$ is called a {\it Borel isomorphism}.

Let $G$ be a group with identity $1$, $\X$ a set. Consider a partially defined function
$m\colon  G\times \X\rightarrow \X,\,\,(g,x)\mt m(g,x)=g\cdot x\in \X.$
For convenience, we shall write $\ex g\cdot x$ to mean that $(g,x)$ belongs to  the domain of $m.$ Then, following \cite{KL}, $m$  is called a (set theoretic) {\it partial action} of $G$ on $\X,$ if for all $g,h\in G$ and $x\in \X$ the following assertions hold:
\smallskip

\noindent (PA1) $\ex g\cdot x$ implies $\ex g\m\cdot (g\cdot x)$ and $g\m\cdot (g\cdot x)=x,$
\smallskip

\noindent (PA2)  $\ex g \cdot (h\cdot x)$ implies $\ex (g h)\cdot x$ and  $g \cdot (h\cdot x)=(g h)\cdot x,$
\smallskip

\noindent (PA3) $\ex 1\cdot x,$ and $1\cdot x=x$.
 \smallskip

We fix the following notations.
\begin{itemize}

\item  $G*\X =\{(g,x)\in G\times \X\mid \ex g\cdot x \},$
the domain of $m.$
\item $\X_{g\m}=\{ x\in \X\mid \ex g\cdot x\},$ and $m_g\colon \X_{g\m}\ni x\mt g\cdot x\in\X_g.$
\item $G^y=\{g\in G\mid \ex g\cdot y\}.$
\item   For $U\subseteq \X,$ the set $G\cdot U=\{g\cdot u\mid g\in G^u,\, u\in U\}$ is called the {\it saturation} of $U.$
\end{itemize}
A partial action $m\colon G*\X\to\X$ induces a family of bijections $\{m_g\colon \X_{g\m}\to \X_g\}_{g\in G},$ and we denote $m=\{m_g\colon \X_{g\m}\to \X_g\}_{g\in G}.$  We remind the next.

\begin{prop} \label{fam}\cite[Lemma 1.2]{QR} A partial action $m$ of $G$ on $\X$ is a family $m=\{m_g\colon \X_{g\m}\to\X_g\}_{g\in G},$ where $\X_g\subseteq \X,$  $m_g\colon\X_{g\m}\to \X_g$ is bijective, for  all $g\in G,$  and such that:
\begin{itemize}
\item[(i)]$\X_1=\X$\,\  and \,\, $m_1=\rm{id}_\X;$
\item[(ii)]  $m_g( \X_{g\m}\cap \X_h)=\X_g\cap \X_{gh};$
\item[(iii)] $m_gm_h\colon \X_{h\m}\cap  \X_{ h\m g\m}\to \X_g\cap \X_{gh},$ and $m_gm_h=m_{gh}$ in $ \X_{h\m}\cap  \X_{g\m h\m};$
\end{itemize}
for all $g,h\in G.$
\end{prop}

The following are two natural examples of partial actions that will be used later in the paper.

\begin{exe}\label{indu} {\bf Induced partial action:}
Let $u \colon G\times \Y\to \Y$ be an action of $G$ on $\Y$ and $\X\subseteq \Y.$ For $g\in G,$ set $\X_g=\X\cap u_g(\X)$ and let $m_g$ be the restriction of $u_g$ to $\X_{g\m}.$  Then $m\colon G* \X\ni (g,x)\mt m_g(x)\in \X $ is a partial action of $G$ on $\X.$  In this case we say that $m$ is induced by $u.$ 
\end{exe}

\begin{exe}\label{hatmh}  Let  $m$ be a partial action of $G$ on $\X,$ then the family $\hat{m}=\{\hat{m}_g\colon (G\times \X)_{g\m}\to (G\times \X)_{g}\}_{g\in G},$ is a partial action of $G$ on $G\times \X,$ where  $(G\times \X)_{g}=G\times \X_{g}$ and $\hat{m}_g(h,x)=(hg\m, g\cdot x),$ for any $g\in G$ and $(h,x)\in  (G\times \X)_{g\m}.$

\end{exe}

From now on $G$ will denote a topological group and $\X$ a topological space. We consider   $ G\times \X$ with the product topology and $G* \X$  with the subspace topology. Moreover,  $m\colon G*\X\to \X$ will denote a   partial action.

\begin{defi} A topological partial action of  $G$ on  $\X$ is a partial action $m=\{m_g\colon \X_{g\m}\to\X_g\}_{g\in G}$ on the underlying
set $\X$, such that each $\X_g$ is open in $\X$, and each $m_g$ is a homeomorphism, for each $g\in G$. If $m\colon G* \X\to \X$ is continuous, we say that the partial action is continuous.
\end{defi}

Now we introduce the concept of a globalization of a partial action and its  enveloping space. Following \cite[(3) Definition 1.1]{AB} we give the next.

\begin{defi}\label{isopa}Let $m=\{m_g\colon \X_{g\m} \to \X_g\}_{g\in G},$  $\te=\{\te_g\colon \Y_{g\m} \to \Y_g\}_{g\in G}$ be topological partial actions of $G$  on the spaces $\X$ and $\Y$ respectively. A morphism $\phi\colon m\to \te$ consists of a continuous map  $\phi\colon \X\to \Y$ such that $\phi(\X_g)\subseteq \Y_g$  and  the square
\begin{equation}\label{cd}\begin{CD}
\X_{g\m} @>m_g>> \X_g\\
@VV\phi V @VV\phi V\\
\Y_{g\m} @>\te_{g}>> \Y_{g}.
\end{CD}\end{equation}
 is commutative, for all $g\in G.$
\end{defi}

The category {\bf TopPar(G)}  has as objects the class of topological partial actions and morphisms described as above.   We say that  $m$ is {\it equivalent} to  $\te$ if  there is an isomorphism between $m$ and $\te.$

A {\it globalization} of a  $m$  is a pair $\{(\beta, \Y),\iota\},$ where $\beta$ is an action  on a space $\Y$ and   $\iota\colon \X\to \Y$ is an injection such that $m$ is equivalent to   $\beta\restriction\iota(\X) .$ A globalization $\{(\beta, \Y),\iota\}$ of $m$ is said to be {\it minimal} if for any globalization  $\{(\beta', \Y'), \iota')\} $ of $m$, there exists a continuous monomorphism $\lambda \colon \beta \to \beta '$ in  ${\bf TopPar(G)}.$

It was shown independently by Abadie \cite[Theorem 1.1]{AB} and  J. Kellendonk and M. Lawson \cite[Theorem 3.4]{KL}, that any topological partial action of $G$ on $\X$ has a minimal globalization, which is  denoted by $\X_G
,$ and called the {\it enveloping space}. We recall their construction.

Let $m$ be a topological partial action of $G$ on $\X.$
Consider  the following
equivalence relation on $ G\times \X.$
\begin{equation}\label{eqgl}(g,x)R  (h,y) \Longleftrightarrow x\in \X_{g\m h}\,\,\,\, \text{and}\, \,\,\, m_{h\m g}(x)=y,\end{equation}
and denote by  $[g,x]$   the equivalence class of the pair $(g,x).$
The  enveloping space  of  $\X$ is the set $\X_G=(G\times \X)/R$  endowed with the quotient topology, so the quotient map
\begin{equation}\label{que}q\colon G\times\X\ni (g,x)\to[g,x]\in \X_G\end{equation} is continuous and open (see \cite[Proposition 3.9 (ii)]{KL}). Now let $\iota \colon \X\ni x\mt [1,x]\in \X_G$.

Then $m$ is equivalent to a partial action induced by  the
action $\mu \colon G\times \X_G\ni (g, [h,x])\to[gh,x]\in  \X_G$  on $\iota(\X).$ The map $\mu$ is  called the {\it enveloping action}.

The following is the basic result about the enveloping space.

\begin{teo}
\label{Abadie} \cite[Theorem 1.1]{AB} and  \cite[Theorem 3.12]{KL}. Let $m$ be a continuous partial action of $G$ on $\X$. Then
\begin{enumerate}
\item  $\iota : \X\rightarrow  \iota(\X)$ is a homeomorphism.
\item $m$ is equivalent to the partial action induced on $\iota(\X)$ by the global action $\mu$.
\end{enumerate}
\end{teo}

\begin{rem} Let $m=\{m_g\colon \X_{g\m} \to \X_g\}_{g\in G},$ and  $\te=\{m_g\colon \Y_{g\m} \to \Y_g\}_{g\in G}$ be objects in {\bf TopPar(G)}. If $m$ is equivalent to $\te,$ then the enveloping spaces $\X_{G}$ and $\Y_{G}$ are homeomorphic. Indeed, let $\phi\colon m\to \te$  be an isomorphism  in  {\bf TopPar(G)}. Then the map    $\phi\colon \X\to \Y$ is an  homeomorphism. Moreover, follows from \eqref{cd}, that  map $\X_G\ni [g,x]\mapsto  [g,\phi(x) ]\in \Y_G$ is well defined and continuous, thus  it is a homeomorphism.
\endproof

\end{rem}

\begin{rem} There is a monomorphism $ m  \to \hat{m}$ in  {\bf TopPar(G)}. Indeed, let $j  \colon \X \ni x\mapsto (1,x) \in G\times \X, $ then $j$ is a continuous injection. Moreover, by the definition of $\hat{m}$ we have that $(G\times \X)_g=G\times \X_g,$  and $j( \X_g )\subseteq G\times \X_g ,$ or any $g\in G$.  It is not difficult to see that diagram \eqref{cd} commutes.
We also remark that in general, $m$ and $\hat{m}$ are not isomorphic. For instance, if $G=\Z$ and $\X$ is the Cantor set, then $G\times \X$ is not  homeomorphic to $\X$ because it is not compact.
\end{rem}

\section{The enveloping space and the orbit equivalence relation}
In this section we show the enveloping space $\X_G$ is homeomorphic to a quotient of $G\times\X$ respect to the orbit equivalence relation given by a partial action of $G$ on $G\times \X$.  We will use that result in the next section to study the metrizability of  $\X_G$.

\begin{defi}  Let  $m$ be  a topological partial action of  $G$ on $\X.$ The orbit equivalence relation $E_G^p$ on $\X$ is given by
$$x E_G^p y \Longleftrightarrow \ex\, g\in G\,\, (g\cdot x\,\,\,\textrm{is defined and}\,\,\, g\cdot x=y).$$
\end{defi}
The set  $\X/E_G^p$  is endowed with the quotient topology. Thus  the quotient map
\begin{equation}\label{pix}\pi_{\X} \colon \X\ni x\mapsto [x] \in \X/E_G^p\end{equation}  is continuous.

\begin{lema}\label{piopen} Let $m=\{m_g\colon \X_{g\m}\to\X_g\}_{g\in G}$ be a   topological partial action of  $G$ on  $\X.$ Then  $\pi_{\X}$ is open.
\end{lema}
\proof Let $U$ be an open set of $\X,$ then
\begin{align*}
\pi_{\X}\m(\pi_{\X}(U))&=\{x\in \X\mid [x]=[u],\,\,\, \text{for some}\,\,\, u\in U\}\\
&=\{x\in \X\mid x=g\cdot u,\,\,\, \text{for some}\,\,\, u\in U \,\,\, \text{and}\,\,\,g\in G^u\}\\
&=\bigcup\limits_{g\in G}m_g(U\cap \X_{g\m}),
\end{align*}
then $\pi_{\X}(U)$ is open in $ \X/E_G^p.$
\endproof

In Example \ref{hatmh} we introduced a partial action $\widehat{m}$ of $G$ on $G\times \X$ associated to the  partial action $m$. We denote by  $\widehat{E}^p_G$ the  orbit equivalence relation induced by $\hat{m}$. The following result shows that $\X_G$ can be described as a quotient space by an orbit equivalence relation.

\begin{teo}\label{mtomg}Let  $m$ be  a topological partial action of a  group $G$ on a space $\X.$ Then  $\X_G$ is the quotient space of $G\times \X$ by $\widehat{E}^p_G.$
\end{teo}

\proof Notice that
\begin{align*} (g, x)R(h, y)& \Longleftrightarrow x \in \X_{ g\m h}\,\, \, \text{and}\,\,\, m_{h\m g}(x) = y \\ &\Longleftrightarrow \exists f \in G\, ( f \cdot x = y \,\, \, \text{and}\,\,\, g f \m= h )\\
  &\Longleftrightarrow \exists f \in G \,( f \cdot  (g,x) = y\,\, \, \text{and}\,\,\, g f \m= h )\\
 &\Longleftrightarrow(g, x) \widehat{E}^p_G(h, y) .
\end{align*}
\endproof

It is known that, in general,  $\X_G$ is not Hausdorff (see for instance \cite [Example 1.4]{AB}). Now we present some results about this problem. The following  is probably known but we give its proof for the sake of completeness.

\begin{lema}\label{closedh} Let $\X$ be a  Hausdorff space and $\rho$ an equivalence relation on $\X$ such that the quotient map $f\colon\X\to \X/\rho $ is open. Then  $\X/\rho$ is Hausdorff, if and only if, $\rho$ is closed in $\X\times \X.$
\end{lema}
\proof Suppose that $\rho$ is closed in $\X\times \X,$ and take $x,y\in \X$  non equivalent points  under $\rho,$ i.e. $(x,y)\notin \rho.$ Since $\rho$ is closed and $\X\times\X$ has the product topology, there exist basic open sets $U,V$ in $\X$ such that $(x,y)\in U\times V$ and $(U\times V)\cap\rho=\emptyset.$ Then  $f(x)\in f(U),$  $f(y) \in f(V)$ and $f(U)\cap f(V)=\emptyset$. Since $f$ is open, then  $\X/\rho$ is Hausdorff.

Conversely, if $u,v\in \X$ are not $\rho$-equivalent,  then $f(x)\neq f(y)$ and there are open sets $\mathcal{O}_u$ and  $\mathcal{O}_v$ in  $\X/\rho$ that separate $f(x)$ and $f(y),$ let $U=f\m(\mathcal{O}_u)$ and $V=f\m(\mathcal{O}_v),$ then  $(u,v)\in U\times V$ and  $(U\times V)\cap\rho=\emptyset,$ which implies that $\rho$ is closed as desired.\endproof

\begin{lema}
\label{Ehat-cerrado}Let  $m$ be  a topological partial action of a  metric topological group $G$ on a metric  space $\X.$
If $G*\X$ is closed, then $\widehat{E}^p_G$ is closed and consequently $\X_G$ is Hausdorff.
\end{lema}

\proof
Let $(g_n, x_n) \widehat{E}^p_G (h_n, y_n)$ such that $g_n\rightarrow g$,  $h_n\rightarrow h$, $x_n\rightarrow x$ and $y_n\rightarrow y$. We will show that  $(g, x) \hat{E}^p_G (h, y)$.
By definition of $\widehat{E}^p_G$, for each $n\in \N$, there is $u_n\in G$ such that  $ g_nu\m_n=h_n$, $x_n\in \X_{u\m_n}$ and  $u_n\cdot x_n=y_n$.  Let $u=h\m g$, then $u_n\rightarrow u$. Since $G*\X$ is closed and $(u_n, x_n) \in G*\X$, one has that $(u,x)\in G*\X$. Therefore  $x\in \X_{u\m}$, $gu\m=h$ and   $u\cdot x=y$. Hence $(g, x) \widehat{E}^p_G (h, y)$ and $\widehat{E}^p_G$ is closed, finally by  Lemmas \ref{piopen} and \ref{closedh}, the space $\X_G$ is Hausdorff.
\endproof

Notice that by the definition of a topological  partial action, each $\X_g$ is open, therefore when $G*\X$ is closed, then they are actually clopen.  This condition was already used in \cite[Proposition 2.1]{EG}  (see also \cite[Proposition 1.2]{AB}), where it was shown that for $\X$ the Cantor space and $G$ a countable discrete group, $\X_G$ is Hausdorff if, and only if, each $\X_g$ is clopen.

\section{Metrizability of $\X_G$}

Now we will address the question of when $\X_G$ is metrizable or Polish. We start  with the following result.

\begin{teo}\label{glopol} Let $m$ be a topological partial action of a separable metrizable  group $G$ on a separable metrizable space $\X.$ Then the following assertions are equivalent.
\begin{enumerate}
\item  $\X_G$ is metrizable.
\item   $\X_G$ is regular and $T_1$.
\end{enumerate}
Moreover, if $G$ and $\X$ are Polish, then any of those conditions are equivalent to $\X_G$ being Polish.

\proof It is clear that ${\rm (1)}\Rightarrow  {\rm (2)}.$   To see that ${\rm (2)} \Rightarrow  {\rm (1)}$,  notice that in this case $\X_G$ is separable and second countable, because  the quotient map $q\colon G\times \X\to \X_G$ is open (\cite[(ii) Proposition 3.9]{KL}). Therefore, by  Urishon's metrization Theorem, $\X_G$ is metrizable.  It remains to prove that if  $\X$ and $G$  are Polish  and $\X_G$ is metrizable, then $\X_G$ is Polish. A classic theorem of Sierpinski  says  that any metrizable space which is the  continuous and open image of a Polish space is also Polish (for a proof see  \cite[Theorem 2.2.9]{GA}).  This is the case of $\X_G$, as the quotient map is open and a product of Polish spaces is also Polish.
\endproof
\end{teo}

By Lemma \ref{Ehat-cerrado} we  know that $\X_G$ is Hausdorff when $G*\X$ is closed. Now we will address the question of the regularity of $\X_G$.
We will use a  generalization of   a theorem of Mackey-Hjorth. We state their result.

\begin{teo}
\label{MH}(Mackey-Hjorth) \cite[Theorem 2.3.5]{KEB} Let $H$ be a Polish group and $G\subseteq H$ a closed subgroup of $H$. Let $\X$ be a Polish $G$-space with an action $a(g,x)$. Then there is a Polish $H$-space $\Y$ with an action $b(h,y)$ such that $\X$ is closed  in $\Y$ and $a(g,x)=b(g,x),$ for $x\in \X$ and $g\in G$.
\end{teo}

Although we only need a particular  case (namely, $H=G$) we present it in a general form to stress  the  connection between  globalizations of partial actions and   extensions of group actions.    First, we  prove some auxiliary results.

\begin{lema}\label{contsep} If $m$ is a  topological partial action,  then $m$ is open.
\end{lema}
\proof Since  $m_g\colon \X_{g\m}\to \X_g$ is a homeomorphism for all $g\in G,$ then
for any open set $U\subseteq \X_{g\m},$ the set
$g\cdot U=m_g(U)$ is also open. Now we check that $m$ is open.
Let $W\subseteq  G*\X$ be  an open set, then there are  families $\{A_i\}_{i\in I}$ and  $\{W_i\}_{i\in I}$  of open sets
in $G$ and $\X$ respectively, such that $W=\bigcup\limits_{i\in I}(G*\X)\cap (A_i\ti W_i).$ For $i\in I$, one has
$$
m((G*\X)\cap (A_i\ti W_i))=\bigcup\limits_{g\in A_i}g\cdot
(W_i\cap \X_{g\m})
$$
which is open in $\X.$ We conclude that $m$ is open.\endproof

\begin{lema}
\label{interiores} Let $m$ be a continuous  partial action  of $G$ on  $\X$ such that $G*\X$ is open.
For each $M\subseteq G$, set  $\X^{M}=\bigcap_{g\in M} \X_{g^{-1}}$.  Then for any  open neighborhood  $N$ of $1$ and  any $x\in \X$,  there is an open neighborhood $M$ of $1$ such that $M\subseteq N$ and  $x\in  int(\X^{M}),$ the interior of  $\X^{M}.$
\end{lema}

\proof Let $N$  be an open neighborhood  of $1$ and  $x\in \X$. Since $G*\X$ is open and $m$ is continuous at $(1,x)$, there are an open set $V$ of $\X$ with $x\in V$ and  an open  neighborhood $M$ of $1$ such that $M\subseteq N$ and    $M\times V\subseteq G*\X$.  Then $(g,  y)\in G*\X$ for all $g\in M$ and $y\in V$.
Therefore $y\in\X_{g^{-1}}$  for all $g\in M$ and $y\in V$. Hence $x\in V\subseteq \X^{M}$ and $x\in int(\X^{M})$.
\endproof

Now we introduce a partial action similar to  $\hat{m}.$
\begin{exe}\label{hatmh2}
Let $H$ be a group,  $G$ be a subgroup of $H$  and  $m$ be a partial action of $G$ on $\X.$ Then the family $\hat{m}^H=\{\hat{m}^H_g\colon (H\times \X)_{g\m}\to (H\times \X)_{g}\}_{g\in G},$ is a partial action of $G$ on $H\times \X,$ where  $(H\times \X)_{g}=H\times \X_{g}$ and $\hat{m}^H_g(h,x)=(hg\m, g\cdot x),$ for any $g\in G$ and $(h,x)\in  (H\times \X)_{g\m}.$ In particular, if $H=G$, then $\hat{m}^H=\hat{m}$.

\end{exe}

 Now we show the main ingredient  for the metrizability of $\X_G$. Recall that by Lemma \ref{piopen} and  Lemma \ref{closedh}, $\X_G$ is Hausdorff iff $\widehat{E}^p_G$ is closed.

\begin{teo}
\label{regcuo}
 Let $G$ and $H$ be  topological groups with $G\subseteq H$  and $\X$ a topological regular space. Suppose that  $m$ is  a continuous  partial action of $G$ on  $\X$ such that  $G*\X$ is open.
Let  $\hat{m}^H$ be the  partial action of $G$ on $H\times \X$ given in Example \ref{hatmh2}.
Then  $(H\times\X)/\widehat{E}^p_H$ is regular under any of the following conditions:
\begin{enumerate}
\item  $\X_g$ is closed, for all $g\in G$.

\item $\X$ is a locally compact metric space and $\widehat{E}^p_G$ is closed in $(G\times\X)^2$.

\item $\X$ is compact Hausdorff and $\widehat{E}^p_G$ is closed in $(G\times\X)^2$.
\end{enumerate}

\end{teo}

\proof

 We first show that (3) implies that every $\X_g$ is closed, for each $g\in G,$ and thus  we  are in case  (1). Let $g\in G$ and $x\in \X_{g}^\mathsf{c}.$ Then for any $y\in \X$ the pair $((1,x),(g\m,y))$ is not in $\widehat{E}^p_G$, and  there are open sets $N_y,M_y$ of $G$ and $V_y,W_y$ of $\X$ such that $(1,x)\in N_y\times V_y,$ $(g\m,y)\in M_y\times W_y$ and $((N_y\times V_y)\times  (M_y\times W_y))\cap \widehat{E}^p_G=\emptyset. $ Since $\X$ is compact, there exist $y_1,\ldots, y_n\in \X$ such that $\X=\bigcup_{i=1}^n W_{y_i}.$ Consider the open set  $V=\bigcap_{i=1}^nV_{y_i}.$ Then $x\in V$ and we shall prove that $V\subseteq \X_{g}^\mathsf{c}. $ Indeed, if   $z\in V\cap \X_{g}$, then there is  $i_0\in \{1,\ldots, n\}$ such that  that $((1,z),(g\m,g\m\cdot z))\in ((N_{y_{i_0}}\times V_{y_{i_0}})\times  (M_{y_{i_0}}\times W_{y_{i_0}}))\cap\widehat{E}^p_G, $ which is a contradiction. From this we conclude that $\X_{g}^\mathsf{c}$ is open, for every $g\in G$.

Now we consider cases (1) and (2).
Let  $F\subseteq H\times \X$  be a  $\hat m^H$-invariant closed set and $(h,x)\notin F.$
It suffices to find a $\hat m^H$-invariant open set  $\mathcal{O}_1$ containing  $(h,x)$  and an open set $\mathcal{O}_2$ containing  $F$ such that $\mathcal{O}_1\cap \mathcal{O}_2=\emptyset.$ Since $F$ is $\hat m^H$-invariant and closed, we can fix $V$ an open subset of $\X$ with $x\in V$ and an open neighborhood  $N\subseteq H$ of the identity  such that
\begin{equation}
\label{v}G\cdot (hN\times V)\cap F=\emptyset.
\end{equation}
By Lemma \ref{interiores},  there is an open  neighborhood of the identity $M\subseteq N$ such that  $x\in  int(\X^{M\cap G})$. So we assume that $N=M$. Since  $G*\X$ is open in $G\times\X$, then $G^x$ is open in $G$. As $1\in G^x$ and $m$ is continuous at $(1,x)$, there are open sets $N_1$ and $V_1$  such that $1\in N_1,$  $x\in V_1,$  and
\begin{equation}
\label{v2}
N_1^2\cap G\subseteq N\cap G^x\,\,\,\,\,\,\,\,\,\,{\rm and }\,\,\,\,\,\,\,\,\,\,\overline{m(N_1^2 \times  V_1\cap G*\X)}\subseteq V.
\end{equation}

Let $N_2\subseteq H$ be  a symmetric open set containing $1$ such that $h\m N_2^2h\subseteq N_1,$ in particular  $h\m N_2h\subseteq N_1.$ Notice also that $\X^{N\cap G}\subseteq \X^{N_1^2\cap G}$ and hence $x\in int(\X^{N_1^2\cap G})$.

Define  $\mathcal{O}_1=G\cdot (N_2h\times (V_1\cap int(\X^{N_1^2\cap G})),$ then  $(h,x)\in \mathcal{O}_1$ and $\mathcal{O}_1$ is an invariant open set (as $\hat m$ is open by  Lemma \ref{contsep}).
Let $(l,y)\in F.$ We need to find an open set  $\mathcal{O}_2$ with $(l,y)\in \mathcal{O}_2$ and $\mathcal{O}_1\cap \mathcal{O}_2=\emptyset.$  We consider two cases:

{\bf Case 1:} $lG\cap hN_1=\emptyset$.  First we claim that  $N_2h\cap N_2lG=\emptyset$. Otherwise, there are $b\in G$, $p,q\in N_2$ such that $plb=qh$. Then $lb\in N_2^2h\subseteq hN_1$ which contradicts our assumption.
Now  we claim that   $\mathcal{O}_2= N_2l\times \X$ works. Indeed, suppose $\mathcal{O}_1\cap \mathcal{O}_2\neq\emptyset$.  Then
there are $v\in V_1\cap int(\X^{N_1^2\cap G})$, $a\in G^v$, $p,q\in N_2$ such that  $pha\m=ql$, which contradicts that $N_2h\cap N_2lG=\emptyset$.

\medskip

{\bf Case 2:} $lG\cap hN_1\not =\emptyset$.  Let $g\in G$ and  $r\in N_1$ such that $lg\m=hr$.

\medskip

{\bf Subcase 2a:}
$g\in G^y .$ Then $g\cdot (l, y) \in F$ and  $g\cdot (l, y)\notin hN \times V,$ thanks to \eqref{v}.  Since $lg\m=hr\in hN,$ one gets that $g\cdot y\notin V.$
Thus, by \eqref{v2}, there exists an open neighborhood $U$ of $g\cdot y$ with $U\cap m(N_1^2 \times  V_1\cap G*\X)=\emptyset.$ The open set  $\mathcal{O}_2=hN_1g \times  m_{g\m} (U\cap \X_g) $ contains $(l,y)$ and  we claim it is disjoint from $\mathcal{O}_1$. If $\mathcal{O}_1\cap \mathcal{O}_2\neq\emptyset$, then there exist $v\in V_1\cap int(\X^{N_1^2\cap G}),$ $u\in U\cap \X_g,\, p\in N_2$ and $a\in G^v$ such that $pha\m\in hN_1g$    and $a\cdot v=g\m \cdot u.$ In particular
$$u=g\cdot (a\cdot v)=(ga)\cdot v=m(ga,v),$$ where the first equality is (PA1) and the second one  is a consequence of (PA2). Since $N_2$ is symmetric, we get that $ga\in N_1h\m N_2h\subseteq N_1^2$. Then,  $u \in U\cap m(N_1^2 \times  V_1\cap G*\X)=\emptyset$, which is a contradiction.

\smallskip

{\bf Subcase 2b} $g\in G\setminus G^y.$  Then $y\in \X_{g\m}^\mathsf{c}$.  We   treat (1) and (2) separately.

\medskip

(1) Suppose  that $\X_{g\m}$ is clopen, so it suffices to show that $\mathcal{O}_2=hN_1g\times \X_{g\m}^\mathsf{c}$ is disjoint from  $\mathcal{O}_1.$
Otherwise, as in the previous case, there are $v\in V_1\cap int(\X^{N_1^2\cap G})$, $a\in G^v$ and $p\in N_2$  such that  $pha\m\in hN_1g,$    and $a\cdot v\notin  \X_{g\m}.$
As before,  we get that $ga\in N_1^2$.  Thus $a\in  g^{-1}N_1^2$. Let $n_1\in N_1^2$ be such that $a=g^{-1}n_1$.  Notice that  $n_1\in G$ as $a,g\in G$ and $G$ is a group. Since  $v\in \X_{N_1^2\cap G}\cap \X _{a^{-1}}\subseteq\X_{n^{-1}_1}\cap \X ^{n^{-1}_1g}$ it follows by Proposition  \ref{fam} (ii) that $n_1\cdot v\in\X_{n_1}\cap \X_{g}.$ Hence $a\cdot v\in \X_{g^{-1}}$, which is a contradiction.

 (2)  Suppose $\X$ is a locally compact metric space and $\widehat{E}^p_G$ is closed.
We claim that there is an open set $W$ such that $y\in W$ and $g\cdot (W\cap \X_{g\m})\cap \overline{V}=\emptyset$. Suppose not and let $\{W_i\}_i$ be an open local base for $y$. Then for each $i$, there is $y_i\in W\cap \X_{g\m}$ such that $g\cdot y_i\in \overline{V}$. As $\overline{V}$ is compact, we assume w.l.o.g. that $g\cdot y_i$ converges to some point $z\in \overline{V}$. Then  $(1,y_i)\hat{E}^p_H (g\m, g\cdot y_i),$ for all $i$. Since $\hat{E}^p_H$ is closed, we get that $(1,y)\hat{E}^p_H (g\m, z)$, therefore $g\in G^y$, which is a contradiction.

Let $W$ be as in the previous claim and  $O_2=h N_1g\times W$. It suffices to show that $O_1$ and $O_2$ are disjoint.
Otherwise, as in  case (1), there are $v\in V_1\cap int(\X^{N_1^2\cap G})$, $a\in G^v$ and $p\in N_2$  such that    $pha\m\in hN_1g,$    and $a\cdot v\in W.$
    As before,  we get that $ga\in N_1^2$.  Take  $n_1\in N_1^2$  such that $a=g^{-1}n_1,$ then    $n_1\in G.$ Since  $v\in \X^{N_1^2\cap G}\cap \X _{a^{-1}}\subseteq\X_{n^{-1}_1}\cap \X _{n^{-1}_1g}$ it follows by Proposition  \ref{fam} (ii) that $n_1\cdot v\in\X_{n_1}\cap \X_{g}.$ Hence $a\cdot v\in \X_{g^{-1}}$ and $g\cdot(a\cdot v)=n_1\cdot v\in V$ by \eqref{v2}. Which contradicts the choice of $W$.  \endproof

\medskip

If in the previous theorem $m$ is a global action, then $(H\times\X)/\hat{E}^p_H$  is the space $\Y$ mentioned in the conclusion of Theorem \ref{MH}. On the other hand, if $H=G$ we get, by Theorem \ref{mtomg}, the enveloping space $\X_G$.

\medskip

 Finally, we have the main result of the paper.

\begin{teo}\label{xgpol}  Let $m$ be a continuous partial action of a separable metrizable group $G$ on a separable metrizable space  $\X$, then $\X_G$ is metrizable under any of the following conditions:
 $G*\X$ is clopen or $\X$ is locally compact and $\widehat{E}^p_G$ is closed. If in addition,  $G$ and $\X$ are Polish, then  $\X_G$ is Polish.

\end{teo}
\proof  By letting $H=G$ in Theorem \ref{regcuo}, then $(G\times\X)/\widehat{E}^p_G$ is  $\X_G$ (by Theorem \ref{mtomg}) and the conclusion follows from Lemma \ref{Ehat-cerrado} and  Theorem \ref{glopol}.
\endproof

 The following example shows that, in general, $\X_G$ is not necessarely Hausdorff even when $G$ is Polish,  $\X$ is a compact Polish space and  $G*\X$ open.

\begin{exe}\label{nonh}  Consider the partial action  of the discrete group  $\Z$ on $\X=[0,1]$ given by $m_{0}=\rm{id}_\X,$ $m_{n}=\rm{id}_V,$ for $n\neq 0,$ and
  $V=(0,1].$  Notice that
$$\Z*\X=\bigcup_{n\in \Z} \{n\}*\X=\{0\}* \X\,\,\cup\bigcup_{n\in \Z\atop n\neq 0} \{n\}*\X=\{0\}\times \X \,\cup\bigcup_{n\in \Z\atop n\neq 0} \{n\}\times V,$$ is open in $\Z\times \X$. Since $V$ is not closed in $\X,$ we have by  \cite[Proposition 1.2]{EG} that $\X_\Z$ is not Hausdorff. In fact, it is easy to verify that the equivalence classes of the form $[(n,0)]$ cannot be separated by open sets.  One may also notice that points in $\X_\Z$ are closed, so $\X_\Z$ is $T_1.$ \end{exe}

\section{ A universal globalization for partial actions}
\label{universal}

Becker and Kechris \cite[Theorems 2.6.1, 5.2.1]{KEB} have shown that for every Polish group $G$ there is Polish space $\UU_G$ and an action of $G$ on $\UU_G$ which is universal, that is, for every Polish space $\X$ and any action of $G$ on $\X$ there is a (necessarily invariant) set $\Y\subseteq \UU_G$ such that $\X$ and $\Y$ are Borel isomorphic as $G$-spaces. In particular, given a  continuous  partial action $m$ of  $G$ over $\X,$ with $G*\X$ clopen in $G\times \X$, since the space  $\X_G$ is Polish (Theorem \ref{xgpol}),  $\X_G$ is Borel isomorphic to a restriction of $\UU_G$.  Therefore there is $\Y\subseteq \UU_G$ (which is now not necessarily invariant) such that $m$ is Borel isomorphic to the $G$ action on $\UU_G$ restricted to $\Y$. In this section we explore those ideas in the context of partial actions of discrete countable groups following the presentation given in \cite{KEB}.

If the function $\phi$ in Definition \ref{isopa} is  bijective and algebraically an isomorphism but is just Borel measurable  rather than  continuous, we will say that the partial actions $m$ and $\theta$ are {\em Borel isomorphic}.

Let $G$ be  a discrete group. By the usual identification of a set with its characteristic function, we may see the collection of all subsets of $G$ as $\{0,1\}^G$ with the product topology.   A natural action of $G$ on $\{0,1\}^G$ is given by
\begin{equation*}   G\times \{0,1\}^G\ni (g,F)\mt g\cdot F=\{gh: h\in F\}\in \{0,1\}^G. \end{equation*}
It is easy to verify such an action is continuous.   Now consider the countable product $(\{0,1\}^G)^\N$ with the product topology and the pointwise product global action $\theta$, more specifically,
$$\theta\colon G\times (\{0,1\}^G)^\N \ni [(g,  (F_n)_n)]\to (g\cdot F_n)_n\in (\{0,1\}^G)^\N.$$  It is easy to see that $\theta$ is a continuous action.
\smallskip

With the notations above we have.

\begin{teo}
\label{product-action}
Let $m$ be a continuous partial action of a countable  discrete group $G$ on a Polish space $\X$ such that $\X_g$ is clopen, for all $g\in G$. Then there is a Borel measurable embedding $\pi: \X\rightarrow (\{0,1\}^{G})^\N$ such that  $m$ and $\theta \restriction \pi[\X]$ are Borel isomorphic, where  $\theta \restriction \pi[\X]$ denotes the induced partial action of $\theta$ on $\pi[\X].$ In particular,  $\X_G$ is Borel isomorphic to $G\cdot \pi[\X]$, the saturation of $\pi[\X]$ inside   $(\{0,1\}^{G})^\N$.
\end{teo}

We need the following

\begin{lema}\label{gxg}Let $m$ be a partial action of a group $G$ on a set $\X.$ Then for any $g\in G$ and $x\in \X_g$ we have $G^xg=G^{g\m\cdot x}.$
\end{lema}
\proof Take $h\in G^x, $ then $\exists h\cdot x$ and  by (PA1) $\exists h\cdot [g\cdot (g\m \cdot x)],$ and (PA2) implies that $\exists (h g)\cdot (g\m \cdot x),$ thus  $G^xg\subseteq G^{g\m\cdot x}.$ Conversely, for $u\in G^{g\m\cdot x}$ one has that $\exists u\cdot (g\m\cdot x)$ and by (PA2)  $\exists (ug\m)\cdot x,$ from this we have that $u=(ug\m)g,$ with $ug\m\in G^x.$
\endproof

\noindent{\em Proof of Theorem \ref{product-action}:}
Let $(V_n)_n$ be a countable base  for $\X$ and consider the following maps $\pi_n: \X\rightarrow \{0,1\}^G$
\[
\pi_n(x)=\{h^{-1}:\; h\in G^{x}\;\;\; \&\;\; \;h\cdot x\in V_n\}.
\]
Let $\pi(x)=(\pi_n(x))_n$. We claim $\pi$ is the required embedding.  To see that each $\pi_n$ is continuous, let  $h\in G$ and $n\in \N$, then
$\{x\in \X:\; h\in \pi_n(x)\}=\{x\in \X_{h^{-1}}: h\cdot x\in V_n\}$ is  open as $h$ is continuous. Similarly,   $\{x\in \X:\; h\not\in \pi_n(x)\}$ is closed. Therefore  $\pi$ is continuous.

To see that  $\pi$ is injective, let $x,y\in\X$, $x\neq y$ and  $n\in\N$ such that $x\in V_n$ and $y\notin V_n$. Then $1\in \pi_n(x)\setminus \pi_n(y)$. Since $\X$ is Polish and $\pi$ is Borel measurable and injective, then   $\pi[\X]$ is Borel  and $\pi$ is a Borel isomorphism between $\X$  and $\pi[X]$ (see \cite[Corollary 15.2]{KE}).

Now we will show that $\pi$ and its inverse are morphisms. Let $g\in G$ and $x \in \X_{g^{-1}}$ we will prove  that $\pi(g\cdot x)=g\cdot \pi(x)$, that is, $\pi_n(g\cdot x)=g\cdot \pi_n(x),$ for all $n\in\N$. By  Lemma \ref{gxg},  $h\in G^{g\cdot x} $
iff $h\in G^x g^{-1}$  iff $h^{-1}=gj^{-1},$ for some $j\in G^x $. Also, $h\cdot (g\cdot x)\in V_n$ iff $j\cdot x\in V_n$.  To show that the inverse of $\pi$ is also a morphism, it suffices to see that given $x,y\in \X$ and $g\in G$ such that $\theta(g,\pi(x))=\pi(y)$, then $(g, x)\in G*\X$. In fact, let $G^x_n=\{h\in G^x:\; h\cdot x\in V_n\}$. Since $\theta(g,\pi(x))=\pi(y)$, then  $g\pi_n(x)=\pi_n(y)$ for all $n\in \N$. That is to say, $g(G^x_n)^{-1}=(G^y_n)^{-1}$ for all $n\in \N$. Thus, $G^x_n =(G^y_n)g$ for each $n\in \N$. It is clear that $G^z=\bigcup_n G^z_n$ for all $z\in \X$.  Therefore,  $G^{x}=G^{y}g$. Since $1\in G^{y}$, then $g\in G^{x}$ and we are done.

Finally, we show that $\X_G$ is Borel isomorphic to $G\cdot \pi[\X]$. Since $G$ has the discrete topology and each $\X_g$ is clopen, then $G*\X$ is clopen. Consider the map $F:X_G\rightarrow G\cdot \pi[\X]$ given by $F([g,x])=g\cdot \pi(x)$. Then $F$ is a continuous bijection. Since by Theorem \ref{xgpol} $\X_G$ is Polish, then as before we conclude that  $F$ is a Borel isomorphism. \qed
\medskip

There is a particular case where we can show that a Polish  space with a continuous partial action is in fact isomorphic to a restriction of $\{0,1\}^G$.

A collection $\{A_j\}_j$ of subsets of $\X$ {\em separates points}, if for all $x,y\in\X$ with $x\neq y$, there is $j$ such that $x\in A_j$ and $y\not\in A_j$ or $x\not\in A_j$ and $y\in A_j$.

\begin{teo}
\label{product-action2}
Let $m$ be a continuous partial action of a countable discrete group $G$ on a Polish  space $\X$ such that $\X_g$ is clopen for all $g\in G$. Suppose that $\{\X_g\}_g$ separates points, then there is a Borel measurable embedding $\pi: \X\rightarrow \{0,1\}^{G}$ such that  $m$ is Borel isomorphic to $\theta_1 \restriction \pi[\X]$, where $\theta_1$ is the global action on  $\{0,1\}^{G}$ given by $\theta_1(g,F)=gF,$  for all $F\in  \{0,1\}^{G}$ and $g\in G$. In particular, $\X_G$ is Borel isomorphic to the saturation of $\pi[\X]$ inside $\{0,1\}^{G}$.
\end{teo}

\proof
Consider the following map $\pi: \X\rightarrow \{0,1\}^G$, $\pi(x)=\{h^{-1}:\; h\in G^{x}\}$. To see that $\pi$ is injective, let $x,y\in\X$, $x\neq y$ and $g\in G$ such that  $x\in \X_g$ and $y\notin \X_g$. Then $g\in \pi (x)\setminus \pi(y)$.  The rest of the proof is analogous to that of  Theorem \ref{product-action}.\endproof

Observe that Theorem \ref{product-action2} is not useful for global actions, however there are interesting examples of partial actions where it can be applied as we show next.

\begin{exe}

 Let $\X=\{0,1\}^\Z$ and $h:\{0,1\}^\Z\rightarrow \{0,1\}^\Z$ be the shift given by $h(x)(n)=x(n+1)$.
Consider the usual global action of $\Z$ on $\X$ given by $h\colon \Z\times \X\ni (n,x)\mapsto h^n(x)\in \X$.  Let $\W=\{x\in \X: \; x(0)=0\}$ and $m$ be the induced partial action on $\W$. That is,
$$
\W_n=\W\cap h^n(\W)=\{x\in \X:\; x(0)=x(-n)=0\},
$$
and $m_n\colon \W_{-n}\ni x\to h^n(x)\in \W_n,$ for $n\in \Z$.  Observe that $\W$ is a clopen subset of $\X$ and thus it is compact. We will show that the family  $\{\W_n\}_{n\in\N}$ separates points. Let $x,y\in\W$ with $x\neq y$ and  $k\in \Z$ be such that $x(k)\neq y(k)$. W.l.o.g we assume that $x(k)=0$ and $y(k)=1$. Then $h^{-k}(x)\in \W$ and $h^{-k}(y)\not\in \W$. Thus $x\in \W_k$ and $y\not\in \W_k$. Let $\pi:\W\rightarrow   \{0,1\}^\Z$ be as in the proof of Theorem \ref{product-action2}. Notice that $\Z ^x=\{n\in \Z:\; x(n)=0\},$ for each $x\in\W,$ and $Y$ be the saturation of $\pi[\W]$ by the global action of $\Z$ on  $\X$, then $Y=\Z \cdot \pi[\W]$. Thus $Y=\{0,1\}^\Z\setminus\{\emptyset\}$. Therefore $\pi[\W]=\{A\in \{0,1\}^{\Z}:\; 0\in A\}$ is an open subset of $Y$ and thus $Y$ is homeomorphic to the globalization $\W_\Z$ by  \cite[(2) Theorem 3.13]{KL}.
\end{exe}

\noindent {\bf Acknowledgment:} We are very thankful to the referee for all his (her)  comments that improved the presentation of the paper and, in particular, for an observation  he (she) made which  clarified Theorem \ref{mtomg}.

\end{document}